\documentclass[11pt]{amsart}
\usepackage{cite}
\usepackage{amssymb,latexsym,amsmath}
\usepackage{amsfonts}

\newtheorem{theorem}{\bf Theorem }[section]

\renewcommand \ker {{\rm ker}\,}

\newcommand \Deg {{\rm Deg}\,}
\renewcommand \dim {{\rm dim}\,}
\newcommand \grad {{\rm grad}\,}
\newcommand \ind {{\rm ind}\,}

\newfont{\wncyr}{wncyr10}  

\newfont{\msbm}{msbm10 at 12pt}
\newcommand \N {\mbox{\msbm \symbol{78}}}
\newcommand \R {\mbox{\msbm \symbol{82}}}

\newcommand \D {\displaystyle}

\newcommand \ff {{\bf f}}
\newcommand \finfty {{\bf f'_{\infty}}}

\newcommand \g {{\bf g}}
\newcommand \ginfty {{\bf g'_{\infty}}}

\newcommand \go {{\bf g'(} 0 {\bf )}}

\newcommand \q {{\bf q}}
\newcommand \qinfty {{\bf q'_{\infty}}}

\newcommand \qo {{\bf q'(} 0 {\bf )}}

\def\cal{\mathcal}

\newcommand \W {\stackrel{\!\!\! \circ}{W _2 ^1}}
\newcommand \Ws {\Bigl[{\stackrel{\!\!\! \circ}{W _2 ^1}}\Bigr]^*}

\newcommand \fxi {f'_{\xi}}

\newcommand \Fuv {\D\int\limits_\Omega
\langle f(x, \gradu), \gradv \rangle_{\R^n} dx}

\newcommand \Guv {\D\int\limits_\Omega g(x,u) vdx}

\newcommand \Puv {\D\int\limits_\Omega
\langle p(x) \gradu, \gradv \rangle_{\R^n} dx}
\newcommand \Quv {\D\int\limits_\Omega
\langle q(x,u), \gradv \rangle_{\R^n} dx}

\newcommand \Lminus {L_{\frac{2n}{n-2}}}

\newcommand \Ln {L_{(n)}}

\renewcommand \div {{\rm div}\,}
\newcommand \gradu {{\rm grad}\, u}
\newcommand \gradv {{\rm grad}\, v}

\newcommand \Hrp {{\cal H}_{r,p}}
\newcommand \Heven {{\cal H}_{r,p}^{even}}
\newcommand \Hodd {{\cal H}_{r,p}^{odd}}
\newcommand \Hevenkn {{\cal H}_{k,n}^{even}}
\newcommand \Hoddkn {{\cal H}_{k,n}^{odd}}

\newcommand \Hevenk {{\cal H}_{k,1}^{even}}
\newcommand \Hoddk {{\cal H}_{k,1}^{odd}}
\newcommand \Hevenl {{\cal H}_{l,1}^{even}}
\newcommand \Hoddl {{\cal H}_{l,1}^{odd}}

\begin{document}
\dedicatory{Dedicated to the Memory of Prof. Igor V. Skrypnik}

\title[The Skrypnik degree]{THE SKRYPNIK DEGREE THEORY AND BOUNDARY VALUE
PROBLEMS}
\author{A.P. Kovalenok \and P.P. Zabreiko}
\address{Institute of Mathematics of NAS of Belarus, 11 Surganov Str.,
220072 Minsk, Belarus}

\begin{abstract}
The paper presents theorems on the calculation of the index of a
singular point and at the infinity of monotone type mappings.
These theorems cover basic cases when the principal linear part of a
mapping is degenerate. Applications of these theorems to proving
solvability and nontrivial solvability of the Dirichlet problem for
ordinary and partial differential equations are considered.
\end{abstract}

\maketitle
\section*{Introduction}

Boundary value problems (BVP) for ordinary and partial
differential equa\-tions constitute one of the most natural
branches where the homotopic theory of monotone type mappings is
successfully applied. It seems to be F. Browder \cite{Browder9}
who suggested a natural construction of reducing a boundary value
problem to some operator equation with a monotone type mapping.
The homotopic theory of monotone type mappings was suggested by
I.V. Skrypnik \cite{S2} (see also \cite{S3}) and developed by many
authors (see for ex. \cite{Browder1,Ber3}). It appeared to be a
natural generalization of the classical Brouwer and Hopf theory of
finite dimensional mappings (see for ex. \cite{GeomMeth}) and
found its numerous application in qualitative studies of BVP (see
for ex. \cite{S3}).

The Krasnosel'skii technique of investigating BVP based on
calculation of topological characteristics of Leray and Schauder
type mappings at zero points and infinity is well known
(see for ex. \cite{GeomMeth,TopMeth}).
But the class of problems where monotone type mappings are applied
is much more broader than that leading to Leray and
Schauder mappings. Moreover, the scheme of reducing
BVP to an operator equation with a monotone type mapping
is absolutely natural and does not require any extra constructions
as such a scheme, concerning Leray and Schauder maps, does.

It is the aim of the paper to develop a technique based on
calculating topological characteristics at zero and infinity
of monotone type mappings for proving the solvability
and existence of nontrivial solution of BVP.
In this connection, in part 1 of the paper, we present theorems for
calculation of the index of zero and infinity of quasimonotone
mappings in Hilbert spaces. In part 2 we consider applications of
theorems stated in part 1 to obtaining conditions for solvability
and existence of nontrivial solutions of the Dirichlet problem
for the 2nd order differential equations.

\section{Topological characteristics of quasimonotone mappings}

\subsection{Preliminaries}

Let $X$ be a Hilbert space, $X^*$ be its dual and
let $\langle l, x \rangle$ denote the pairing
between the elements $x\in X$ and $l\in X^*$.

A mapping $\Phi$ defined on a set $D\subseteq X$
with values in the space $X^*$ is said to be

{\it demicontinuous} if it maps each sequence $x_n \in D$
converging to $x_0$ into a sequence $\Phi x_n$ weakly converging
to $\Phi x_0$ ($\Phi x_n \stackrel{w}{\to} \Phi x_0$);

{\it quasimonotone} if for each sequence $x_n \in D$ with
$x_n \stackrel{w}{\to} x_0$
the inequality
$$\limsup_{n \to \infty} \ \langle \Phi x_n, x_n - x_0\rangle \ge 0$$
holds;

{\it of class $(S_+)$} or to have {\it $(S_+)$ property}
if for each sequence $x_n \in D$ with
$x_n \stackrel{w}{\to} x_0$
the inequality
$$\limsup_{n \to \infty} \ \langle \Phi x_n, x_n - x_0\rangle \le 0$$
implies that $x_n \to x_0$.

The class of mappings with $(S_+)$ property
defined on a weakly closed set $D$ involves
{\it strictly monotone} mappings, i.e. such that
$$\langle \Phi x -\Phi y, x - y \rangle \ge m ||x-y||_X \ \ \
(x,\  y \in D)$$
for some  $m>0$. Quasimonotone perturbations of $(S_+)$ mappings
have $(S_+)$ property as well. Furthermore, any
compact mapping is quasimonotone. So the class of mappings
having $(S_+)$ property and
defined on a weakly closed set $D$ contains
compact perturbations of strictly monotone mapping.

For any triple $(\Phi,D,y)$, where $y\in X^*$ and $\Phi$ is a
demicontinuous quasimo\-notone mapping defined on the closure
$\overline{D}$ of the open bounded set $D\subseteq X$ and
satisfying the condition
$$\inf\limits_{x \in \partial\Omega}||\Phi x - y ||>0, $$
one associates an integer topological characteristic ---
{\it the Skrypnik degree $\Deg (\Phi,D,y)$
of the mapping $\Phi$ of the set
$D$ with respect to the point $y$}.
As well as the classical Brouwer and Hopf and Leray and Schauder
degrees the Skrypnik degree is determined by its properties:
normality with respect to $J^{-1}$
(we let $J$ denote a Riesz isomorphism between $X$ and $X^*$),
additivity under domain and homotopy invariance (see for ex. \cite{S3}).
Moreover, in Hilbert spaces the Skrypnik degree theory generalizes
the classical Leray and Schauder degree theory (see for ex. \cite{Ber3}).

A point $x_0$ is a {\it zero of a mapping $\Phi$} if $\Phi x_0 =0$.
A zero $x_0$ is called {\it isolated}
if there exists a ball $B_{r_0}(x_0)$ which does not contain
any other zeroes of the mapping $\Phi$.
A zero $x_0$ is called {\it strictly isolated} if
there exists $r_0 >0$ such that $\Phi x_n \to 0$
and $||x_n - x_0|| \le r_0$ imply $x_n \to x_0$.

Let $x_0$ be a strictly isolated zero of a demicontinuous
quasimonotone mapping $\Phi$ defined in some neighbourhood of
$x_0$. Then the degree \linebreak $\Deg (\Phi,B_r(x_0),0)$ is the
same for all balls $B_r(x_0)$ with  the center $x_0$ and a
sufficiently small radius $r$, this common degree being referred
as {\it the index of the zero $x_0$} and denoted by $\ind
(x_0,\Phi)$. Let us note that to define the index $\ind
(x_0,\Phi)$ for a demicontinuous mapping with $(S_+)$ property it
suffices for $x_0$ to be an isolated zero.

If for a mapping $\Phi$ defined on all elements of the space $X$
with large norms there exists $r_0 > 0$ such that
$\Phi x_n \to 0$ and $||x_n||\ge r_0$ imply
$x_n \to \infty$ then $\Phi$ is said to be {\em strictly
nondegenerate at infinity}. For any demicontinuous and quasimonotone
mapping $\Phi$ which is strictly nondegenerate at infinity
the degree $\Deg (\Phi,B_r(x_0),0)$ is the same for all
balls $B_r(x_0)$ with an arbitrary center $x_0$ and a sufficiently large
radius $r\ge r_0$, this common degree being referred as
{\it the index at infinity} or {\it asymptotical index}
of the mapping $\Phi$ and denoted by $\ind (\infty,\Phi)$.
It is clear that to define the index at infinity for a demicontinuous
mapping with $(S_+)$ property it suffices for the set of zeroes
of this mapping to be bounded.

Formulas for calculation of the index play the very important part
when studying operator equations. For mappings of monotone type
the Kronecker theorem holds (see for ex. \cite{S3}). The
corollary of this theorem is the fact that (under corresponding assumptions)
the asymptotical index of a mapping equals the sum of zeroes indexes
of this mapping. So index formulas make one able to
prove the solvability and existence of nontrivial solutions,
bifurcation points etc.

\subsection{Index of a zero}

Let $\Phi$ be a demicontinuous and quasimonotone map\-ping defined
in a neighbourhood of its zero point $x_0$. We assume $\Phi$ to be
Fr{\'e}chet differentiable  at $x_0$, i.e.
$$\Phi (x_0 + h)=\Phi'(x_0)h + \omega(h), \ \ \
\lim\limits_{h\to 0}\frac{\omega(h)}{||h||_{X^*}}=0,$$
with the Fr{\'e}chet derivative $\Phi'(x_0)$ being
quasimonotone and such that the mapping $J\Phi'(x_0)$
proves to be a Fredholm operator.

The space $X$ is the direct sum of the subspaces $X_1$ and $X^1$
invariant for $J\Phi'(x_0)$ where $X_1$ is the finite dimensional
root subspace of the operator $J\Phi'(x_0)$ and $X^1$ is complement to $X_1$.
Let $X_0 = \ker \Phi'(x_0)$ and $X^0$ be complement to
$X^0$. It is evident that $X_0 \subseteq X_1$
and without loss of generality we can suppose $X^1\subseteq X^0$.
Let $P_0$, $P^0=I-P_0$ and $P_1$, $P^1=I-P_1$ denote
projectors onto $X_0$, $X^0$ and $X_1$, $X^1$ respectively.
At last, note that there exists a nondegenerate linear
operator $T$ which acts in $X_1$ and such that
$TJ\Phi'(x_0)P_1=P^0P_1$ (see for ex. \cite{ZK2}).

Now we are in position to formulate the following

\begin{theorem}\label{DegZero} {\rm \cite{ZabKova1,ZabKova2,Preprint}}
Let the operator $\omega(h)$ have a
positively homogeneous of order $l>1$ principal term at zero, i.e.
$$\omega(h)= C_m(h) +\omega_m(h),$$
with
$C_m$ being a positively homogeneous of order  $l>1$
continuous mapping and
$||\omega_m(h)||_{X^*}=o(||h||_X^m)$
when $||h||_X \to 0$.

If $0$ is the unique zero point
of the finite dimensional mapping $\Theta_m u= P_0 TP_1 J C_m u$
($u\in X_0$) then $x_0$ is a strictly isolated
zero of the mapping $\Phi$ and
$$
\ind (x_0, \Phi) =
(-1)^{n_0-l} \cdot \ind (\Phi'(x_0); X^1)
\cdot \ind (0, \Theta_m; X_0),
$$
when $n_0=\dim X_1$ and $l$ is the number of Jordan blocks in
the Jordan form of $J\Phi'(x_0)$ in the space $X_1$.

If, besides,  $J\Phi'(x_0)$ has only
isolated Fredholm points of the spectrum on
$(-\infty,0)$ then
$$
\ind (x_0, \Phi)
= (-1)^{\nu(J\Phi'(x_0))+n_0 - l} \cdot \ind (0, \Theta_m; X_0),
$$
where $\nu(J\Phi'(x_0))$ is the sum of multiplicities
of negative eigenvalues of $J\Phi'(x_0)$.
\end{theorem}

It should be noted that in papers \cite{ZabKova1,ZabKova2,Preprint}
cited above more general cases were considered.

\subsection{Asymptotical index}

Let $\Phi$ be a demicontinuous and quasimonotone mapping
defined on all elements of $X$ with sufficiently large norms.
We assume $\Phi$ to be asymptotically differentiable
(differentiable at infinity), i.e.
$$\Phi h=\Phi'(\infty)h + \omega(h), \ \ \
\lim\limits_{h\to +\infty}\frac{\omega(h)}{||h||_{X^*}}=0,$$
with the asymptotical derivative $\Phi'(\infty)$ being
quasimonotone and such that the mapping $J\Phi'(\infty)$
proves to be a Fredholm operator.

Similarly to the case of zero point the
space $X$ is the direct sum of the subspaces $X_1$ and $X^1$ where
invariant for $J\Phi'(\infty)$
the $X_1$ is the finite dimensional
root subspace of the operator $J\Phi'(\infty)$ and $X^1$ is complement to $X_1$.
Furthermore, let $X_0 = \ker \Phi'(x_0)$ and $X^0$ be complement to $X^0$
such that $X^1\subseteq X^0$, $P_0$, $P^0=I-P_0$ and $P_1$, $P^1=I-P_1$
denote projectors onto $X_0$, $X^0$ and $X_1$, $X^1$ respectively,
$T$ be nondegenerate operator acting in $X_1$ be such that
$TJ\Phi'(\infty)P_1=P^0P_1$.

\begin{theorem}\label{DegInf} {\rm \cite{Preprint,ZabKova3}}
Let the operator $\omega(x)$ have a
positively homogeneous of order $l>1$ principal term at infinity, i.e.
$$\omega(x)= C_k(x) +\omega_k(x),$$
with $C_k$ being a positively homogeneous of
order $0\le k < 1$ continuous mapping and
$||\omega_k(x)||_{X^*}=o(||x||_X^k)$ when $||x||_X \to +\infty$.

If $0$ is the unique zero point of the finite dimensional mapping
$\Theta_k u= P_0 TP_1 J C_k u$ ($u\in X_0$) then
the mapping $\Phi$  is strictly nondegenerate at infinity
and
$$
\ind (\infty, \Phi) =
(-1)^{n_0-l} \cdot \ind (\Phi'(\infty); X^1)
\cdot \ind (0, \Theta_k; X_0),
$$
when $n_0=\dim X_1$ and $l$ is the number of Jordan blocks in
the Jordan form of $J\Phi'(\infty)$ in the space $X_1$.

If, besides,  $J\Phi'(\infty)$ has only
isolated Fredholm points of the spectrum on
$(-\infty,0)$ then
$$
\ind (\infty, \Phi)
= (-1)^{\nu(J\Phi'(\infty)) + n_0 -l} \cdot
\ind (0, \Theta_k; X_0),
$$
where $\nu(J\Phi'(\infty))$ is the sum of multiplicities
of negative eigenvalues of $J\Phi'(\infty)$.
\end{theorem}

\section{Applications to boundary value problems}

Here we consider the Dirichlet problem (in the weak sense) of the form
\begin{equation}\label{D}
-\div (f(x,\gradu)+q(x,u)) + g(x,u)=0,
\ u\Big|_{x\in \partial\Omega}=0,
\end{equation}
where $f(x,\xi): \Omega \times  \R^n \mapsto \R^n$,
$q(x,t):\Omega \times \R \mapsto \R^n$,
$g(x,t):\Omega \times \R \mapsto \R$,
$\Omega \subset \R^n$.

Let us describe in short the scheme (which is quite standard)
of our further study.
Under natural conditions including restrictions on the growth
of vector functions $f(x,\xi)$, $q(x,t)$ and the function $g(x,t)$ the
problem (\ref{D}) is equivalent to the operator equation
\begin{equation}\label{OpEq}
\Phi u =0
\end{equation}
with the monotone type mapping $\Phi = F + Q + G$ between the
Sobolev space $\W$ and its dual space $\Ws$
where operators $F$, $Q$ and $G$ are given by
\begin{equation}\label{Definition}
\begin{array}{ccc}
\langle Fu,v\rangle =\Fuv,
\\[12pt]
\langle Qu,v\rangle =\Quv, \ \ \
\langle Gu,v\rangle =\Guv
\end{array}
\end{equation}
when $u,\ v \in \W$.
For the problem (\ref{D}) to have a solution (nontrivial solution)
it is sufficient then, due to the Kronecker theorem (see \cite{S3}),
to prove that the asymptotical index of $\Phi$ does not
equal zero (the index of trivial solution).

We present two types of conditions under which
the problem (\ref{D}) is solvable. They differ by
the way of calculating the asymptotical index of $\Phi$.
The first condition is based on the calculation of the index by
the theorem \ref{DegInf} and leads to the equation with the sublinear with
respect to the variable $t$ function $g(x,t)$.
The second condition is the consequence of
a priori estimates on $q(x,t)$ and $g(x,t)$ under which
$\ind (\infty,\Phi)=1$ and allows to consider
equations with the vector function $q(x,t)$ being asymptotically zero
and the function $g(x,t)$ having
a one-side superlinear growth with respect to $t$.

Also we use both these ways of calculating the asymptotical index
to obtain conditions under which the problem (\ref{D}) has
a nontrivial solution. Besides that here
we apply the theorem \ref{DegZero} and the Skrypnik theorem \cite{S3}
to calculate the zero index of $\Phi$.
The Skrypnik theorem
makes one possible to calculate the index when $\Phi$
has just the G\'{a}teaux derivative
(nondegenerate in some neighbourhood of the zero and
satisfying some other conditions).
The theorem \ref{DegZero}
can be applied provided $\Phi$ is the Fr{\'e}chet differentiable.
It is well known that there is no
nonlinear superposition operator acting in $L_2$
which is Fr{\'e}chet differentiable.
That is why possible applications of the theorem \ref{DegZero}
are restricted
by the class of problems with the linear principal part $F$.
However, corresponding results based on this theorem
remain valid for quasilinear equations. So we distinguish for our further
discussions the quasilinear problem
\begin{equation}\label{D'}
-\div (p(x)\gradu+ q(x,u)) + g(x,u)=0,
\ u\Big|_{x\in \partial\Omega}=0
\end{equation}
where $p(x):\R^n \mapsto \R^n$ is a linear
with respect to $\xi$ for any $x\in \Omega$ operator (matrix),
which is the special case of the problem (\ref{D}).

We let $\ff$, $\q$ and $\g$ respectively denote superposition
operators generated by the vector functions $f(x,\xi)$, $q(x,t)$
and the function $g(x,t)$. Further, such properties of
superposition operators as the action between appropriate spaces,
differentiability, asymptotical differentiability, and so on are
needed. All these properties can be expressed in terms of growth
restrictions on $f(x,\xi)$, $q(x,t)$, $g(x,t)$ and their
derivatives (see for ex. \cite{AppellZabreiko}). These growth
restrictions are determined by spaces where the superposition
operator acts and, in our case, depend significantly on embeddings
of the Sobolev space $\W$ into spaces of integrable functions. To
provide the maximal growth, i.e. to cover as wide range  of
problems involved as possible, we assume for the set $\Omega$ to
have such a smooth boundary (see for ex. \cite{FunctionSpaces})
that embeddings
\begin{equation}\label{Embedding2}
\W\subset \Ln,
\end{equation}
where
$$
\Ln = \begin{cases}
C, & \textrm{when}\ n=1 \cr
E_M, & \textrm{when}\ n=2\cr
\Lminus, & \textrm{when}\ n>2,
\end{cases} \  M(s)=e^{s^2}-1,
$$
hold.
Let us note that embeddings (\ref{Embedding2}) are limit
and except the case $n=1$ noncompact.

In what follows we use notation $\Hrp\ (r \in \R, \  p\in \N)$ for
the family of functions $h(x,\eta)$ which are
positively homogeneous of order $r$ with respect to $\eta \in \R^p$.
Furthermore, we let $\Heven,\
\Hodd$ denote respectively the set of
even and odd with respect to $\eta$  functions
$h(x,\eta)$ belonging to $\Hrp$.

\subsection{Reduction}

We consider the problem (\ref{D}) provided that

$A_1)$ there exists $m>0$ such that
$$
\langle f(x,\xi) - f(x,\eta), \xi - \eta  \rangle_{\R^n} \ge
m ||\xi - \eta||_{\R^n}^2 \ \ \ (\xi,\ \eta \in \R^n);
$$

$A_2)$ the superposition operator $\ff$ acts from $L_2^n :=
\underbrace{L_2\times \ldots \times L_2}_{n\ times}$ into $L_2^n$;

$A_3)$ the superposition operator $\q$ acts as the improving one from
$\Ln$ into $L_2^n$;

$A_4)$ the superposition operator $\g$ acts as the improving one from
$\Ln$ into $\Ln^*$.

Due to $A_1)$ --- $A_4)$ the problem (\ref{D})
is equivalent to the operator equation (\ref{OpEq})
with the continuous operator $\Phi$ having $(S_+)$ property
as the sum of the continuous and strongly monotone operator
$F$ and completely continuous operators $Q$ and $G$.

In the case of the linear principal part $F=P$ where
$$
\langle Pu,v\rangle =\Puv
$$
the assumptions
$A_1)$ and $A_2)$ can be written in the form

$A_1)'$  all the elements of the matrix
$p(x)$ belong to $L_\infty$ and there exists $m>0$
such that
$$
\langle p(x)\xi,\xi \rangle_{\R^n}\ge m ||\xi||^2_{\R^n}\ \ \
(x\in \Omega, \xi\in \R^n).
$$

So conditions $A_1)'$, $A_3)$, $A_4)$ are supposed
to be fulfilled when studying the problem
(\ref{D'}).

\subsection{Solvability}

\subsubsection{}
 We first consider the case when the theorem
\ref{DegInf} is applicable. For the mapping $\Phi$ to be
asymptotically differentiable we suppose that for any $x\in \Omega
$ the vector functions $f(x,\xi)$, $q(x,t)$ and the function
$g(x,t)$ have asymptotical derivatives $f'(x,\infty)$,
$q'(x,\infty)$ and $g'(x,\infty)$  which generate asymptotical
deri\-vatives $\finfty$, $\qinfty$ and $\ginfty$ of the
superposition operators $\ff$, $\q$ and $\g$.

The simplest situation when the problem (\ref{D})
has a solution is rather evident and takes place
if the linear problem
\begin{equation}\label{Dinf}
-\div (f'(x,\infty)\gradu+q'(x,\infty)u) + g'(x,\infty)u=0,
\ u\Big|_{x\in \partial\Omega}=0
\end{equation}
has no nontrivial solution, which corresponds to the
case of nondegenerate derivative $\Phi'(\infty)$.
If the resonance phenomenon happens one has to consider
higher order terms at infinity.
The following result
is the consequence of the theorem \ref{DegInf}.

\begin{theorem}\label{solv1}
{\rm \cite{Kova3,Kova4,ZabKova7}} Let for some $0<k<1$ there exist
vector functions $f^k(x,\xi) \in (\Hoddkn)^n:= \underbrace{\Hoddkn
\times \ldots \times \Hoddkn}_{n \ times}$, $q^k(x,t) \in
(\Hoddk)^n$ and a function $g^k(x,t) \in \Hoddk$ which generate
respectively principal terms of order $k$ of the superposition
operators $\ff - \finfty$, $\q - \qinfty$ and $\g - \ginfty$ at
infinity.

Then the problem {\rm (\ref{D})} has a solution
provided that none of nontrivial
solutions of the problem {\rm (\ref{Dinf})}
satisfies the problem
\begin{equation}\label{Dkinf1}
-\div (f^k(x,\gradu)+q^k(x,u)) + g^k(x,u)=0,
\ u\Big|_{x\in \partial\Omega}=0.
\end{equation}
\end{theorem}

\subsubsection{}
 In our next assertion we drop the demand of the asymptotical
differen\-tiability of superposition operators which automatically
implies the sublinea\-rity of their generating functions (see for
ex. \cite{AppellZabreiko}). But the price of this is "the
smallness" of $q(x,t)$ and "the hard" estimate for $g(x,t)$.

\begin{theorem}\label{solv2}
{\rm \cite{ZabKova6}}
Let the superposition operator $\q$ be asymptotically
zero and for some $\D 0< \delta < \frac{m}{K^2}$,
where $K$ denotes the norm of the operator embedding $\W$ into $L_2$,
the estimate
$$
g(x,t)t\ge -\delta t^2\ \ \ (x\in \Omega, \ t\in \R)
$$
hold.

Then the problem {\rm (\ref{D}) } has a solution.
\end{theorem}

\subsection{Existence of nontrivial solutions}

In what follows we assume that the zero of the space $\W$
is the solution of the problem (\ref{D}) and (\ref{D'}).

\subsubsection{}
 Assume that conditions of the subsection 2.2.1 on the
asymptotical differentiability of $\Phi$ are fulfilled. We
consider the case of resonance, when the problem {\rm
(\ref{Dinf})} has nontrivial solutions.

The following result is the consequence of the
theorem \ref{DegInf} and the Skrypnik theorem on the zero index
of the G\'{a}teaux differentiable mapping (see \cite{S3}).

\begin{theorem}\label{nontr1}
{\rm \cite{Kova3,Kova4,ZabKova7}}
Let the following conditions hold

i) the superposition operators $\ff$, $\q$ and $\g$ are
G\'{a}teaux differentiable on $L_2^n$, $\Ln$ and $\Ln$
respectively, with $\fxi (x,0)$, $q'_t(x,0)$ and $g'_t(x,0)$ being
derivatives of the vector functions $f(x,\xi)$, $q(x,t)$  and the
function $g(x,t)$ with respect to the second variable at zero.

ii) for some $0<k<1$ there exist vector functions $f^k(x,\xi) \in
(\Hevenkn)^n$, $q^k(x,t) \in (\Hevenk)^n$ and a function $g^k(x,t)
\in \Hevenk$ which generate respectively principal terms of order
$k$ of the superposition operators $\ff - \finfty$, $\q - \qinfty$
and $\g - \ginfty$ at infinity.

Then the problem {\rm (\ref{D})} has a nontrivial solution
provided that the problem
$$
- \div (\fxi (x,0)\gradu +q'_t(x,0)u) + g'_t(x,0)u= 0,
\ u\Big|_{x\in \partial\Omega}=0
$$
has no nontrivial solutions and
none of nontrivial solutions of the problem {\rm (\ref{Dinf})}
satisfies the problem {\rm (\ref{Dkinf1})}.
\end{theorem}

\subsubsection{}
 Assume that conditions of the subsection 2.2.1 concerning the
asymp\-totical differentiability of $\q$ and $\g$ hold. We also
suppose the superposition operators $\q$ and $\g$ to have
respectively derivatives $\qo$ and $\go$ at zero, with $q'(x,0)$
and $g'(x,0)$ being derivatives of the vector function $q(x,t)$
and the function $g(x,t)$ with respect to $t$ at zero.

We consider here the most sophisticated "resonance" situation
when both the problem
\begin{equation}\label{Dinf2}
-\div (p(x)\gradu+q'(x,\infty)u) + g'(x,\infty)u=0,
\ u\Big|_{x\in \partial\Omega}=0
\end{equation}
and the problem
\begin{equation}\label{Dzero}
- \div (p(x)\gradu +q'(x,0)u) + g'(x,0)u= 0,
\ u\Big|_{x\in \partial\Omega}=0
\end{equation}
has nontrivial solutions, i.e. derivatives of $\Phi$
at zero and infinity are degenerate.

The following assertion is the consequence of the theorem \ref{DegZero}
and the theorem \ref{DegInf}.
It combines two cases, when higher order terms of
$\Phi$ are odd at infinity and even at zero and vice versa.

\begin{theorem}\label{nontr2}
{\rm \cite{Kova3,Kova4,ZabKova7}}
Let the following conditions hold

i) for some $0<k<1$ there exists
a vector function $q^k(x,t) \in (\Hoddk)^n$
($q^k(x,t) \in (\Hevenk)^n$)
and a function $g^k(x,t) \in \Hoddk$
($g^k(x,t) \in \Hevenk$)
which generate respectively principal terms of order $k$ of
the superposition operators $\q - \qinfty$ and $\g - \ginfty$ at infinity;

ii) for some $l>1$ ($\D l < \frac{n}{n-2}$ when $n>2$) there exists
a vector function $q^l(x,t) \in (\Hevenl)^n$
($q^l(x,t) \in (\Hoddl)^n$)
and a function $g^l(x,t) \in \Hevenl$
($g^l(x,t) \in \Hoddl$)
which generate respectively principal terms of order $l$ of
the superposition operators $\q - \qo$ and $\g - \go$ at zero.

Then the problem {\rm (\ref{D'})} has a nontrivial solution
provided that none of nontrivial solutions
of the problem {\rm (\ref{Dinf2})}
satisfies the problem
$$
-\div q^k(x,u) + g^k(x,u)=0,
\ u\Big|_{x\in \partial\Omega}=0
$$
and none of nontrivial solutions of the problem
{\rm (\ref{Dzero})} satisfies the problem
\begin{equation}\label{Dlzero}
- \div q^l(x,u) + g^l(x,u)= 0,
\ u\Big|_{x\in \partial\Omega}=0.
\end{equation}
\end{theorem}

Let us note that in the case when
the problem (\ref{Dzero}) (the problem (\ref{Dinf2}))
turns out to have no nontrivial solutions we have to
consider even higher order terms at infinity (at zero) only.
It seems to be quite clear how to formulate the same results
as the theorem above for such cases (see \cite{Kova3,Kova4,ZabKova7}).
That is why we omit the precise formulating here.

\subsubsection{}
 Similarly to the subsection 2.2.2 we study here the case of the
problem (\ref{D'}) whose function $g(x,t)$ can have one side
superlinear growth with respect to $t$. We assume conditions of
the theorem \ref{solv2} to be fulfilled. Furthermore, we suppose
assumptions of the previous subsection 2.3.2 concerning the
differentiability of $\q$ and $\g$ at zero to hold.

We consider the degenerate case, when the
linearization at zero (\ref{Dzero}) of the problem (\ref{D'})
has nontrivial solutions.

\begin{theorem}
{\rm \cite{ZabKova6}}
Let for some $l>1$ ($\D l < \frac{n}{n-2}$ when $n>2$) there exists
a vector function $q^l(x,t) \in (\Hevenl)^n$
and a function $g^l(x,t) \in \Hevenl$
which generate respectively principal terms of order $l$ of
the superposition operators $\q - \qo$ and $\g - \go$ at zero.

Then the problem {\rm (\ref{D'})} has a nontrivial solution
provided that no solution of the problem
{\rm (\ref{Dzero})} satisfies the problem
{\rm (\ref{Dlzero})}.
\end{theorem}

\subsubsection{}
 As the reader could notice in this section we studied only
situations when the mapping $\Phi$ had the degenerate derivative
at zero or at infinity (or both). We should say that the technique
of proving the existence of nontrivial solutions when both
derivatives (in the subsections 2.3.1 and 2.3.2) or the derivative
at zero (in the subsection 2.3.3) of $\Phi$ are nondegenerate is
quite different. In such cases we do not have to consider higher
order terms but to establish the the parity of the sums of
negative eigenvalues multiplicities of corresponding
linearizations.

So for the existence of a nontrivial solutions of the
problem (\ref{D'}) under assumptions of the subsection 2.3.2
when both the problem (\ref{Dinf2}) and the problem
(\ref{Dzero}) has no nontrivial solutions
it is sufficient to prove that
the sums of negative eigenvalues multiplicities of problems
$$
- \div (p(x)\gradu +q'(x,0)u) + g'(x,0)u= -\lambda \div \grad u,
\ u\Big|_{x\in \partial\Omega}=0,
$$
$$
- \div (p(x)\gradu +q'(x,\infty)u) + g'(x,\infty)u=
-\lambda \div \grad u,
\ u\Big|_{x\in \partial\Omega}=0
$$
have different parities.

In a similar manner one can study analogous situations
if they appear in the subsection 2.3.1 and 2.3.3 (see
\cite{Kova3,Kova4,ZabKova6,ZabKova7}).



\end{document}